\documentclass[a4paper,12pt,reqno]{amsart}
\usepackage[body={17cm,21cm}]{geometry}
\usepackage[initials]{amsrefs}

\usepackage{amssymb,amscd}
\usepackage[all]{xy}

\usepackage[mathscr]{eucal}
\usepackage{enumerate}

\numberwithin{equation}{section}
\theoremstyle{plain}
\newtheorem{thm}[equation]{Theorem}
\newtheorem{cor}[equation]{Corollary}
\newtheorem{lem}[equation]{Lemma}
\newtheorem{prop}[equation]{Proposition}

\newtheorem{definition}[equation]{Definition}
\newtheorem{exa}[equation]{Example}
\newtheorem{rem}[equation]{Remark}

\theoremstyle{definition}

\theoremstyle{remark}

\DeclareMathOperator{\Br}{Br}

\DeclareMathOperator{\Div}{Div}

\DeclareMathOperator{\Hom}{Hom}
\DeclareMathOperator{\id}{id}

\DeclareMathOperator{\Pic}{Pic}

\DeclareMathOperator{\Spec}{Spec}

\def\Br{{\rm Br}}
\def\inv{{\rm inv}}

\def\G{{\mathbb G}}

\def\ker {{\rm  Ker}}

\def\Pic{{\rm Pic}}
\def\ker{{\rm ker}}

\newcommand{\et}{{\textup{\'et}}} 

\DeclareFontEncoding{OT2}{}{} 
   \newcommand{\textcyr}[1]{%
     {\fontencoding{OT2}\fontfamily{wncyr}\fontseries{m}\fontshape{n}%
      \selectfont #1}}

\newcommand{\Sha}{{\mbox{\textcyr{Sh}}}}

\DeclareFontFamily{U}{wncy}{}
\DeclareFontShape{U}{wncy}{m}{n}{%
<5>wncyr5%
<6>wncyr6%
<7>wncyr7%
<8>wncyr8%
<9>wncyr9%
<10>wncyr10%
<11>wncyr10%
<12>wncyr6%
<14>wncyr7%
<17>wncyr8%
<20>wncyr10%
<25>wncyr10}{}
\DeclareMathAlphabet{\cyr}{U}{wncy}{m}{n}
\begin{document}

\title[Descent and Brauer-Manin obstructions]
{Comparing descent obstruction and Brauer-Manin obstruction for open varieties}

\author{Yang CAO}

\address{Yang CAO \newline Laboratoire de Math\'ematiques d'Orsay
\newline Univ. Paris-Sud, CNRS, Univ. Paris-Saclay \newline 91405 Orsay, France}

\email{yang.cao@math.u-psud.fr}

\author{Cyril DEMARCHE}

\address{Cyril DEMARCHE \newline Sorbonne Universit\'es, UPMC Univ Paris 06 \newline Institut de Math\'ematiques de Jussieu-Paris Rive Gauche
\newline UMR 7586, CNRS, Univ Paris Diderot
\newline Sorbonne Paris Cit\'e, F-75005, Paris, France  and
\newline D\'epartement de math\'ematiques et applications
\newline \'Ecole normale sup\'erieure
\newline 45 rue d'Ulm, 75230 Paris Cedex 05, France}

\email{cyril.demarche@imj-prg.fr}

\author{Fei XU}

\address{Fei XU \newline School of Mathematical Sciences, \newline Capital Normal University,
\newline 105 Xisanhuanbeilu, \newline 100048 Beijing, China}

\email{xufei@math.ac.cn}


\date{\today.}


\maketitle

\begin{abstract} We provide a relation between Brauer-Manin obstruction and descent obstruction for torsors over not necessarily proper varieties under a connected linear algebraic group or a group of multiplicative type. Such a relation is also refined for torsors under a torus.  The equivalence between descent obstruction and \'etale Brauer-Manin obstruction for smooth projective varieties is extended to smooth quasi-projective varieties, which provides the perspective to study integral points.
\end{abstract}

\section{Introduction}\label{intro}

The descent theory for tori was first established by Colliot-Th\'el\`ene and Sansuc in \cite{CTS87} and was extended by Skorobogatov to groups of multiplicative type in \cite{Sk}. In a series of papers \cite{H}, \cite{HSk}, \cite{HSk1}, Harari and Skorobogatov introduced descent obstruction for a general algebraic group and compared the descent obstruction with the Brauer-Manin obstruction. By various works of Poonen \cite{P},  the second named author \cite{D09}, Stoll  \cite{St} and Skorobogatov \cite{Sk1},  it was proved that the descent obstruction is equivalent to the \'etale Brauer-Manin obstruction for smooth projective geometrically integral varieties. In this paper, we study the relation between the descent obstruction and the Brauer-Manin obstruction for open varieties by using new arithmetic tools developed in \cite{BD}, \cite{CT08}, \cite{CTX}, \cite{D}, \cite{Ha08} and \cite{HS05}, and we extend the equivalence between the descent obstruction and the \'etale Brauer-Manin obstruction to smooth \emph{quasi-projective} varieties. 

Let $k$ be a number field, $\Omega_k$ the set of all primes of $k$ and ${\bf A}_k$ the adelic ring of $k$. A variety over $k$ is defined to be a separated scheme $X$ of finite type over $k$. Fix an algebraic closure $\bar k$ of $k$. We denote by $X_{\bar k}$ the fibre product $X\times_k \bar{k}$. Let $$\Br(X)=H^2_{\textup{\'et}}(X, \Bbb G_m), \ \ \  \Br_1(X)= \ker(\Br(X) \rightarrow \Br(X_{\bar k}))  \ \ \  \text{and} \ \ \ \Br_0(X)= \rm{Im} (\Br(k) \xrightarrow{\pi^*} \Br(X))  $$ where $X\xrightarrow{\pi} Spec(k)$ is the structure morphism, and $\Br_a(X)=\Br_1(X)/\Br_0(X)$.  For any subgroup $B$ of $\Br(X)$, one can define the Brauer-Manin set
$$ X({\bf A}_k)^B= \{ (x_v)_{v\in \Omega_k}\in X({\bf A}_k): \ \sum_{v\in \Omega_k} \inv_v(\xi(x_v))=0 \ \ \text{for all} \ \xi\in B \}  $$ with respect to $B$. When $B=\Br(X)$, we simply write this Brauer-Manin set as $ X({\bf A}_k)^{\Br}$.

Suppose $Y\xrightarrow{f} X$ is a left torsor under a linear algebraic group $G$ over $k$. The descent obstruction (see \cite{H}, \cite{HSk} and \cite{HSk1}) given by $f$ is defined by the following set
$$ X({\bf A}_k)^f = \{(x_v)\in X({\bf A}_k): ([Y](x_v))\in {\rm{Im}} (H^1(k, G) \rightarrow \prod_{v\in \Omega_k} H^1(k_v, G)) \} = \bigcup_{\sigma\in H^1(k, G)} f_\sigma (Y^\sigma ({\bf A}_k)) $$ where $Y^\sigma \xrightarrow{f_\sigma}X$ is the twist of $Y\xrightarrow{f} X$ by a 1-cocycle representing $\sigma\in H^1(k,G)$. Moreover, one can define
$$ X({\bf A}_k)^{\text{desc}} =\bigcap_{Y\xrightarrow{f} X} X({\bf A}_k)^f$$
following \cite{P}, where $Y\xrightarrow{f} X$ runs through all torsors under all linear algebraic groups over $k$.

The main results in this paper are the following theorems. 

\begin{thm}\label{c-a} {\rm (Theorem \ref{main-general}) } Let $k$ be a number field, $G$ a connected linear algebraic group or a group of multiplicative type over $k$, and $X$ a smooth and geometrically integral variety over $k$. Suppose $Y\xrightarrow{f} X$ is a left torsor under $G$. For any subgroup $A \subseteq \Br(X)$ which contains the kernel of the natural map $f^*: \Br(X) \rightarrow \Br(Y)$ we have
$$ X({\bf A}_k)^A = \bigcup_{\sigma\in H^1(k, G)} f_\sigma (Y^\sigma ({\bf A}_k)^{f_\sigma^*(A)})   $$ where $Y^\sigma \xrightarrow{f_\sigma}X$ is the twist of $Y\xrightarrow{f} X$ by $\sigma$ and $ \Br(X) \xrightarrow{f_\sigma^*} \Br(Y^\sigma)$ is the associated pull-back map, for each $\sigma\in H^1(k,G)$.
\end{thm}

When $G$ is a torus, this theorem can be refined in order to get Theorem \ref{tor} in \S 4. In particular, we prove:

\begin{thm}\label{intor} {\rm (Corollary \ref{algebraic}) } Under the same assumptions as in Theorem \ref{c-a}, if $G$ is assumed to be a torus, then
$$ X({\bf A}_k)^{\Br_1(X)} = \bigcup_{\sigma\in H^1(k, G)} f_\sigma (Y^\sigma ({\bf A}_k)^{\Br_1(Y^\sigma)}) $$
and 
$$ X({\bf A}_k)^{\Br} = \bigcup_{\sigma\in H^1(k, G)} f_\sigma (Y^\sigma ({\bf A}_k)^{\Br_1(Y^\sigma)+f_\sigma^*(\Br(X))})  . $$ 
\end{thm}

This result is inspired by some lectures by Yonatan Harpaz. It should be pointed out that the first part in Theorem \ref{intor} was first obtained by Dasheng Wei in \cite{Wei}: his proof uses an argument of Harari and Skorobogatov in \cite{HaSk} together with an exact sequence due to Sansuc (see \cite{BD}, Theorem 2.8). Theorem \ref{intor} can be applied to study strong approximation, as in \cite{Wei}. It should be noted that in general, the image of $\Br (X)$ in $\Br(Y^\sigma)$ in Theorem \ref{c-a} and Theorem \ref{intor} is not easy to describe, even under the assumption $\bar k[X]^\times=\bar k^\times$ (see \cite[Theorem 1.7(b)]{HSk03}). 

\begin{definition} Let $X$ be a variety over a number field $k$ and let $B$ be a subgroup of $\Br(X)$. For a finite subset $S$ of $\Omega_k$, we denote by $pr^S: X({\bf A}_k) \to X({\bf A}_k^S)$ the projection map,  where ${\bf A}_k^S$ is the set of adeles of $k$ without $S$-components.

We say that  $X$ satisfies \emph{strong approximation off $S$} if $X({\bf A}_{k})\neq \emptyset$ and the diagonal image of $X(k)$ is dense in $pr^{S}(X({\bf A}_{k}))$.

We say that  $X$ satisfies \emph{strong approximation with respect to $B$ off $S$} if $X({\bf A}_{k})^{B} \neq \emptyset$ and the diagonal image of $X(k)$ is dense in $pr^{S}(X({\bf A}_{k})^{B})$.
\end{definition}

Corollary 3.20 in \cite{D} provides a sufficient condition for strong approximation with Brauer-Manin obstruction to hold for a connected linear algebraic group. As an application of Theorem \ref{intor}, we prove that this sufficient condition is also a necessary condition:

\begin{thm} \label{application} {\rm (Corollary \ref{iff})}  Let $G$ be a connected linear algebraic group over a number field $k$ and let $S$ be a finite subset of $\Omega_k$ containing $\infty_k$. Then $G$ satisfies strong approximation with respect to $\Br_1(G)$ off $S$ if and only if $\prod_{v\in S} G'(k_v)$ is not compact for any non-trivial simple factor $G'$ of the semi-simple part $G^{ss}$ of $G$. 
\end{thm}

For any variety $X$ over a number field $k$, one can define, following \cite{P}:
$$ X({\bf A}_k)^{\textup{\'et}, \Br}= \bigcap_{Y\xrightarrow{f} X} \bigcup_{\sigma\in H^1(k, F)} f_\sigma (Y^\sigma({\bf A}_k)^{\Br}) \, ,$$ where $Y\xrightarrow{f} X$ runs through all torsors under all finite group schemes $F$ over $k$. The last two sections of the paper are devoted to the proof of the following generalization of \cite{D09} and \cite{Sk1}:

\begin{thm}\label{inteq} {\rm (Corollary \ref{oneside} and Theorem \ref{main-last}) } If $X$ is a smooth quasi-projective and geometrically integral variety over a number field $k$, then
$$   X({\bf A}_k)^{\rm{desc}} = X({\bf A}_k)^{\textup{\'et}, \Br} . $$
\end{thm}

Terminology and notations are standard if not explained. For any connected linear algebraic group $G$ over an field $k$ of characteristic zero, the reductive part $G^{\rm{red}}$ of $G$ is defined by the exact sequence
$$ 1\rightarrow R_u(G)\rightarrow G \rightarrow G^{\rm{red}} \rightarrow 1 $$ where $R_u(G)$ is the unipotent radical of $G$.  The semi-simple part $G^{ss}$ of $G$ is defined to be the derived subgroup $[G^{\rm{red}}, G^{\rm{red}}]$, which is isogenous to the product of its simple factors, and the maximal toric quotient $G^{tor}$ of $G$ is defined to be $G^{\rm{red}}/[G^{\rm{red}}, G^{\rm{red}}]$.  We use $\hat{G}$ for the character group of $G$. For a topological abelian group $A$, the topological dual of $A$ is defined as $A^D=\Hom_{cont} (A, \Bbb Q/\Bbb Z)$ with the compact-open topology. For any ring $R$, $R^\times$ stands for the group of invertible elements of $R$. For a number field $k$, we denote by $\infty_k$ the set of all archimedean primes of $k$ and by $O_S$ the ring of $S$-integers, for any finite subset $S \subset \Omega_k$ containing $\infty_k$. For any $v\in \Omega_k$, $k_v$ is the completion of $k$ with respect to $v$, and if $v\in \Omega_k\setminus \infty_k$, $O_v$ is the integral ring of $k_v$.

The paper is organized as follows. In \S \ref{sl}, we establish some algebraic results over an arbitrary field of characteristic zero which we need in the next sections. Then we prove Theorem \ref{c-a} in \S \ref{clag}, Theorem \ref{intor} in \S \ref{rtc}. As an application of those results, we prove Theorem \ref{application} in \S \ref{aa}.  Theorem \ref{inteq} is proved in \S \ref{cI} and \S \ref{CII}.

\section{Brauer groups of torsors}\label{sl}

In this section, we assume that $k$ is an arbitrary field of characteristic 0.  
\begin{lem} \label{br} Let $H$ be a semi-simple simply connected group or a unipotent group over $k$. Suppose $X$ is a smooth and geometrically integral variety over $k$.  If $Z\xrightarrow{\rho} X$ is a torsor under $H$, then the induced map $\Br(X)\xrightarrow{\rho^*}\Br(Z)$ is an isomorphism.
\end{lem}

\begin{proof} We first show that $\Br(X)\xrightarrow{\cong} \Br(X\times_k H)$, where the map is induced by the natural projection $X \times_k H \to X$. Using the spectral sequence $$ H^p(k, H^q(X_{\bar k}, \Bbb G_m)) \Rightarrow H^{p+q}(X, \Bbb G_m) , $$ one only needs to show that
$$ \bar k [X_{\bar k}]^\times/\bar k^\times\xrightarrow{\cong} {\bar k}[X_{\bar k} \times_{\bar k} H_{\bar k}]^\times /\bar k^\times , \ \ \ \Pic(X_{\bar k})\xrightarrow{\cong} \Pic(X_{\bar k} \times_{\bar k} H_{\bar k}) \ \ \ \text{and} \ \ \ \Br(X_{\bar k})\xrightarrow{\cong} \Br(X_{\bar k} \times_{\bar k} H_{\bar k}). $$
Since $\bar k[H]^\times=\bar k^\times$ and $\Pic(H_{\bar k})=\Br(H_{\bar k})=0$ by \cite[Proposition 2.6]{CTX}, the first two parts are true by \cite[Proposition 6.10 ]{Sansuc}.  To prove the last part, Kummer exact sequence ensures that one only needs to prove that
\begin{equation} \label{pr}  H^2_{\textup{\'et}}(X_{\bar k}, \Bbb Z/n) \xrightarrow{\cong} H_{\textup{\'et}}^2(X_{\bar k} \times_{\bar k} H_{\bar k}, \Bbb Z/n ) \end{equation} for all $n \geq 1$. This last isomorphism follows from \cite[Proposition 2.2]{SZ} and \cite[Expos\'e XI, Th\'eor\`eme 4.4]{SGA4} with 
$H^i_{\textup{\'et}}(H_{\bar k}, \Bbb Z/n)=0$ for $i=1, 2$. So we proved the required isomorphism $\Br(X)\xrightarrow{\cong} \Br(X\times_k H)$.

Let us now deduce Lemma \ref{br}: since $\Pic(H)=0$, \cite[Proposition 2.4]{BD} gives the following short exact sequence
$$ 0 \rightarrow \Br (X) \rightarrow \Br (Z) \xrightarrow{m^*-p_Z^*} \Br(H\times_{k} Z) \, ,$$ 
where $m^*$ and $p_Z^*$ are induced by the multiplication map $H\times_k Z\xrightarrow{m} Z$ and the projection map $H\times_k Z\xrightarrow{p_Z} Z$ respectively. Since $m\circ (1_{H}\times \id)=p_Z \circ (1_H\times \id ) = \id $, one concludes that $m^*=p_Z^*$ by the above argument. Therefore $\Br(X)\xrightarrow{\cong}\Br(Z)$.
\end{proof}

Let $H$ be a closed subgroup of an algebraic group $G$ over $k$, and $Y\xrightarrow{f} X$ be a left torsor under $H$. Let $Z\xrightarrow{\rho} X$ be the left torsor under $G$ defined by the contracted product $Z=G\times^H Y$ (see \cite[Example 3 in p.21]{Sko}): the torsor $Z$ is the push-forward of $Y$ by the homomorphism $H \to G$. The projection map $G\times_k Y \xrightarrow{pr_G} G$ induces the following commutative diagram
\begin{equation} \label{d} \begin{CD}
G\times_k Y @>>> Z=G\times^H Y  \\
@V{pr_G}VV @VV{\theta}V \\
G @>{\pi}>> G/H \, , \end{CD} \end{equation}
where $\theta$ is induced by $pr_G$ via the quotient by $H$.

\begin{lem} \label{q}  With the above notations, for any $\gamma\in (G/H)(k)$, the composite map  $\theta^{-1}(\gamma) \to Z \xrightarrow{\rho} X$ is naturally a left torsor under $H^\sigma$, which is canonically isomorphic to the twist of $Y\xrightarrow{f} X$ by the $k$-torsor $\pi^{-1}(\gamma)$ under $H$.
\end{lem}
\begin{proof} It follows from diagram (\ref{d}) and \cite[Example 2 in p.20]{Sko}. \end{proof}

Let $G$ be a connected linear algebraic group over $k$, and $Y$ be a smooth variety over $k$. Since $G_{\bar k}$ is rational over $\bar k$ by Bruhat decomposition, the projections $G\times_k Y\to G$ and $G\times_k Y\to Y$ induce an isomorphism 
 $$ \Br_a(G)\oplus \Br_a(Y) \xrightarrow{\sim} \Br_a(G\times_k Y) $$ 
by \cite[Lemma 6.6]{Sansuc}. If $P$ is a (left) torsor under $G$ over $k$ and $H^3(k, \bar{k}^\times)=0$, the previous result generalizes to an isomorphism
\begin{equation} \label{iso-sansuc}
\Br_a(P) \oplus \Br_a(Y) \xrightarrow{\sim} \Br_a(P \times Y)
\end{equation}
by \cite[Lemma 5.1]{BvH}.  

Let $G$ be a connected linear algebraic group over $k$ and let $X$ be a smooth variety over $k$ with $H^3(k, \bar{k}^\times)=0$. Suppose that $Y\xrightarrow{f} X$ is a left torsor under $G$ and $P$ is a left $k$-torsor under $G$, associated to a cocycle $\sigma \in Z^1(k,G)$. One can consider $P$ as a right torsor under $G$ by defining a right action $ x\circ g:= g^{-1} x$ (see \cite[Example 2 in p.20]{Sko}). This right torsor is called the inverse right torsor of $P$ under $G$, and is denoted by $P'$. 
One can now consider the map given by the quotient of $P \times_k Y$ by the diagonal action of $G$ given by $g \cdot (p,y) := (p \circ g^{-1}, g \cdot y) = (g \cdot p, g \cdot y)$:
$$ \chi_P:  P\times_k Y\rightarrow Y^\sigma:=P'\times^G Y \, .$$


\begin{definition} \label{brauer-twist}  With the above notation, assuming that $H^3(k, \bar{k}^\times) = 0$, consider the map
$$ \psi_\sigma = \psi_P : \Br_a(Y^\sigma) \xrightarrow{\chi_P^*} \Br_a (P\times_k Y) \xleftarrow{\sim} \Br_a(P) \oplus \Br_a(Y) \rightarrow \Br_a(Y) \, .$$
\end{definition}

The following lemma, which compares the algebraic Brauer groups of twists of a given torsor, can be regarded as an extension of \cite[Lemma 1.3]{Wei} to torsors under connected linear algebraic groups.

\begin{lem} \label{twist-isom}  The morphism $\psi_\sigma$ in Definition \ref{brauer-twist} is an isomorphism.
\end{lem}
\begin{proof}  
The natural morphism $(pr_P,\chi_P) : P \times_k Y \to P \times_k Y^\sigma$ is an isomorphism, and we have a commutative diagram:  
$$ \begin{CD} 
P\times_k Y @>{(pr_P,\chi_P)}>>  P\times_k Y^\sigma \\
@V{pr_P}VV @VV{pr_P}V \\
P @>>{\id}> P \, .
\end{CD}$$

Therefore $(pr_P,\chi_P)^* : \Br_a(Y^\sigma \times_k P) \to \Br_a(Y \times P)$ induces the identity map on the subgroups $\Br_a(P) \subset \Br_1(Y^\sigma \times_k P)$ and $\Br_a(P) \subset \Br_1(Y \times_k P)$, hence 
$$\psi_\sigma : \Br_a(Y^\sigma) \to \Br_a(Y^\sigma \times_k P) \xrightarrow{(pr_P,\chi_P)^*} \Br_a(Y \times P) \to \Br_a(Y)$$
is an isomorphism (using the isomorphism \eqref{iso-sansuc}).
\end{proof}

Let $f: Y\rightarrow X$ be a torsor under a connected linear algebraic group $G$ over $k$ and let $$ a_Y: \ G \times_k Y\rightarrow Y$$ be the action of $G$. There is a canonical map $\lambda : \Br_1(Y)\rightarrow \Br_a(G)$ by \cite[Lemma 6.4]{Sansuc}. Let $e: \Br_a(G)\rightarrow \Br_1(G)$ be the section of $\Br_1(G)\rightarrow \Br_a(G)$ such that $1_G^*\circ e=0$. If $X$ is smooth and geometrically integral, then the following diagram 
\begin{equation} \label{diag torsor} 
\begin{CD} 
\Br_1(Y)  @>{\lambda}>> \Br_a(G) \\
@VVV @VV{p_G^* \circ e}V \\
\Br(Y) @>>{a_Y^*-p_Y^*}> \Br(G\times_k Y)
\end{CD}
\end{equation}
commutes by \cite[Theorem 2.8]{BD}, where  
$G\times_k Y\xrightarrow{p_G} G$ and $G\times_k Y\xrightarrow{p_Y} Y$ are the projections.  One can reformulate the commutative diagram \eqref{diag torsor} in the following proposition:

\begin{prop}\label{torsorbraueraction}  With the above notation, one has $$ b(t \cdot x)=\lambda (b)(t)+b(x)$$ for any $x\in Y(k)$, $t\in G(k)$ and $b\in \Br_1(Y)$.
\end{prop}

\begin{proof} The commutativity of diagram \eqref{diag torsor} implies that 
$$a_Y^*-p_Y^*=p_G^*\circ e\circ \lambda:  \ \Br_1(Y)\rightarrow \Br_1(G\times Y) \, ,$$
therefore one has 
$$b(t \cdot x)=a_Y^*(b)(t,x)=p_Y^*(b)(t,x)+p_G^*\circ e\circ \lambda (b) (t,x)=b(x)+\lambda (b) (t)$$ as required. 
\end{proof}

\section{Connected linear algebraic groups or groups of multiplicative type}\label{clag}

In this section, we study the relation between the descent obstruction and the Brauer-Manin obstruction for a general connected linear group or a group of multiplicative type. 

First we need the following fact concerning topological groups:

\begin{lem}\label{top} Let $f: M\rightarrow N$ be an open homomorphism of topological groups. If $K$ is a closed subgroup of $M$ containing $\ker(f)$, then $f(K)$ is a closed subgroup of $N$.
\end{lem}
\begin{proof} Since $K$ is a closed subgroup containing $\ker(f)$, one has $$f(K)=f(M)\setminus f(M\setminus K) . $$  Since $f$ is an open homomorphism, $f(M)$ is an open subgroup of $N$. This implies that $f(M)$ is closed in $N$. Since $f(M\setminus K)$ is open in $N$, one concludes that $f(K)$ is closed in $N$.
\end{proof}

\begin{rem} The assumption $K \supseteq \ker(f)$ in Lemma \ref{top} can not be removed. For example, the projection map $pr^S: {\bf A}_k \to {\bf A}_k^S$ is open where ${\bf A}_k^S$ is the set of adeles of $k$ without $S$-component. It is clear that $k$ is a discrete subgroup of ${\bf A}_k$ by the product formula. However $k$ is dense in ${\bf A}_k^S$ 
by strong approximation for $\Bbb G_a$, when $S$ is not empty. 
\end{rem}

For a short exact sequence of connected linear algebraic groups, one has the following result.

\begin{prop}\label{connected-groups} Let $$1\rightarrow G_1 \xrightarrow{\psi} G_2 \xrightarrow{\phi}  G_3 \rightarrow 1$$ be a short exact sequence of connected linear algebraic groups over a number field $k$. Then 

(1)  $\phi\left(G_2({\bf A}_k)^{\Br_1(G_2)}\right)$ is a closed subgroup of $G_3({\bf A}_k)$. 

(2)  If $G'(k_\infty)$ is not compact for each simple factor $G'$ of the semi-simple part of $G_3$, then one has
$$ G_3({\bf A}_k)^{\Br_1(G_3)} = G_3(k) \cdot \phi\left(G_2({\bf A}_k)^{\Br_1(G_2)}\right) . $$ 
\end{prop}

\begin{proof}  Let $S$ be a sufficiently large finite set of primes of $\Omega_k$ containing $\infty_k$ and let
${\bf G}_1$ (resp. ${\bf G}_2$, resp. ${\bf G}_3$) be a smooth group scheme model of $G_1$ (resp. $G_2$, resp. $G_3$) over $O_S$ with connected fibres, such that 
 the short exact sequence of smooth group schemes 
$$ 1 \rightarrow {\bf G}_1 \xrightarrow{\psi} {\bf G}_2 \xrightarrow{\phi} {\bf G}_3 \rightarrow 1$$ extends the given short exact sequence of their generic fibres. The set $H_{\rm et}^1(O_v, {\bf G}_1)$ is trivial by Hensel's lemma together with Lang's theorem, and the following diagram 
$$ \begin{CD} {\bf G}_3(O_v) @>{\partial_v}>> H_{\rm et}^1(O_v, {\bf G}_3) \\
@VVV @VVV \\
G_3(k_v) @>>{\partial_v}> H^1(k_v, G_3) \end{CD} $$ commutes, hence we deduce the following commutative diagram of exact sequences in Galois cohomology:
$$\begin{CD}
G_1(k) @>{\psi}>> G_2(k)  @>{\phi}>> G_3(k) @>{\partial}>> H^1(k, G_1)  \\
@VVV @VVV @VVV @VVV \\
G_1({\bf A}_k) @>>{(\psi_v)}> G_2({\bf A}_k) @>>{(\phi_v)}> G_3({\bf A}_k) @>>{(\partial_v)}>  \bigoplus_{v\in \Omega_k} H^1(k_v, G_1) \, .  \end{CD} $$

In addition, \cite[Theorem 5.1]{D} and \cite[Corollary 6.11]{Sansuc} gives the following commutative diagram of exact sequences of topological groups and pointed topological spaces:
\begin{equation} \label{diag adelic}
\begin{CD}
@. @. G_1({\bf A}_k) @>{\theta_1}>> \Br_a(G_1)^D @>>> \Sha^1(k, G_1) \\
@. @. @VV{(\psi_v)}V @VV{(\psi^*)^D}V @.\\
1 @>>> \ker(\theta_2) @>>> G_2({\bf A}_k) @>{\theta_2}>> \Br_a(G_2)^D   \\
@. @VVV @VV{(\phi_v)}V  @VV{(\phi^*)^D}V \\
1 @>>> \ker(\theta_3) @>>> G_3({\bf A}_k) @>{\theta_3}>> \Br_a(G_3)^D \\
@. @. @VV(\partial_v)V  @. \\
@. @.  \bigoplus_{v\in \Omega_k} H^1(k_v, G_1)  \, , \end{CD}
\end{equation}
where $\Br_a(G_i)^D$ is the topological dual of the discrete group $\Br_a(G_i)$, for $1\leq i\leq 3$. Since $\theta_1(G_1({\bf A}_k))$ is the kernel of the continuous map $\Br_a(G_1)^D \to \Sha^1(k, G_1)$, it is a closed subgroup of $\Br_1(G)^D$. Since $(\psi^*)^D$ is a closed map, one obtains that $(\psi^*)^D(\theta_1(G_1({\bf A}_k))$ is a closed subgroup of $\Br_1(G_2)^D$. It implies that
$$ \ker(\theta_2) \cdot \psi (G_1({\bf A}_k)) = \theta_2^{-1}\left[(\psi^*)^D(\theta_1(G_1({\bf A}_k))\right] $$ is a closed subgroup of $G_2({\bf A}_k)$ by diagram \eqref{diag adelic}. Proposition 6.5 in Chapter 6 of \cite{PR} ensures that $\phi: G_2({\bf A}_k)\rightarrow G_3({\bf A}_k)$ is an open homomorphism of topological groups. Then $\phi(\ker(\theta_2))= \phi\left(G_2({\bf A}_k)^{\Br_1(G_2)}\right)$ is closed by Lemma \ref{top}, and property (1) follows.

Let us now prove statement (2): Corollary 3.20 in \cite{D} (see also the proof of Proposition 4.5 in \cite{CX2}) implies that 
$$\ker(\theta_3)=G_3({\bf A}_k)^{\Br_1(G_3)}=\overline{G_3(k) \cdot G_3(k_\infty)^0}\, ,$$ 
where $G_3(k_\infty)^0$ is the connected component of identity with respect to the topology of $k_\infty$. One only needs to show that
$$ G_3({\bf A}_k)^{\Br_1(G_3)} \subseteq  G_3(k) \cdot \phi\left(G_2({\bf A}_k)^{\Br_1(G_2)}\right) \, . $$
For any $(x_v) \in  \overline{G_3(k) \cdot G_3(k_\infty)^0}$, there is $h\in G_3(k)$ and $h_\infty \in G_3(k_\infty)$ such that $$(\partial_v)(h\cdot h_\infty)= (\partial_v)(x_v) \, ,$$ 
because $(\partial_v)$ is a continuous map with respect to the discrete topology of $\bigoplus_{v\in \Omega_k} H^1(k_v, G_1)$. Since $\phi_\infty (G_2(k_\infty)^0)$ is open and connected, the finiteness of $H^1 (k_\infty, G_1)$ gives
$$G_3(k_\infty)^0=\phi_\infty (G_2(k_\infty)^0) \, .$$ 
Therefore
$$(h\cdot h_\infty) \in G_3(k)  \cdot \phi\left(G_2({\bf A}_k)^{\Br_1(G_2)}\right)  $$ and one can replace $(x_v)$ by $(h\cdot h_{\infty})^{-1} \cdot (x_v)$. Without loss of generality,
one can therefore assume $(\partial_v)(x_v)$ is the trivial element in $\bigoplus_{v\in \Omega_k} H^1(k_v, G_1)$.

Since $\Sha^1 (k, G_1)$ is finite, one can fix $\xi_1, \cdots, \xi_n$ in $G_3(k)$ such that each element of $\Sha^1(k,G_1) \cap \partial (G_3(k))$ is represented by one of the $\xi_i$'s.
As $\partial_\infty (h_\infty)$ is trivial for any $h_\infty \in G_3(k_\infty)^0$, one concludes that
$$(x_v) \in \overline{\bigcup_{i=1}^n \xi_i\phi(\ker(\theta_2))}= \bigcup_{i=1}^n \xi_i \cdot \overline{\phi(\ker(\theta_2))} \subseteq  G_3(k) \cdot \phi\left(G_2({\bf A}_k)^{\Br_1(G_2)}\right) $$ by Corollary 1 in Page 50 of \cite{Ser} and assertion (1). 
\end{proof}

The main result of this section is the following theorem:

\begin{thm}\label{main-general} Let $X$ be a smooth and geometrically integral variety and let $G$ be a connected linear algebraic group or a group of multiplicative type over a number field $k$. Suppose that $f: Y\rightarrow X$ is a left torsor under $G$. If $A$ is a subgroup of $\Br(X)$ which contains the kernel of the natural map $f^*: \Br(X) \rightarrow \Br(Y)$,  then 
$$ X({\bf A}_k)^A = \bigcup_{\sigma\in H^1(k, G)} f_\sigma \left(Y^\sigma ({\bf A}_k)^{ f_\sigma^*(A)}\right) \, ,$$  
where $Y^\sigma \xrightarrow{f_\sigma}X$ is the twist of $f$ by $\sigma$ and $ \Br(X)\xrightarrow{f_\sigma^*} \Br(Y^\sigma)$ is the associated pull-back morphism, for each $\sigma\in H^1(k,G)$. 
\end{thm}

\begin{proof} By the functoriality of Brauer-Manin pairing, one only needs to show that 
$$ X({\bf A}_k)^{A}\subseteq \bigcup_{\sigma\in H^1(k, G)} f_\sigma \left(Y^\sigma ({\bf A}_k)^{ f_\sigma^*(A)}\right) \, .$$
It is clear that
 \begin{equation}\label{equiv}
 (x_v)\in \bigcup_{\sigma\in H^1(k, G)} f_\sigma (Y^\sigma ({\bf A}_k))  \ \ \ \Leftrightarrow \ \ \
 ([Y](x_v))\in {\rm {Im}} \left[H^1(k,G) \rightarrow \prod_{v\in \Omega_k} H^1(k_v, G) \right] \, . \end{equation}

(1) Assume that $G$ is connected.

Recall first that Hensel's lemma together with Lang's theorem ensures that $H^1(k, G)$ maps to $\bigoplus_{v\in \Omega_k} H^1(k_v, G)$. Since any element $P \in \Pic (G)$ can be given the structure of a central extension of algebraic groups
 \begin{equation} \label{ses}  1 \to \Bbb G_m \to P \to G \to 1 \end{equation}  by \cite[Corollary 5.7]{CT08}, one obtains a coboundary map
 $$ \partial_P:  \ \ \ H^1(X, G) \rightarrow H^2(X, \Bbb G_m)=\Br(X)$$ 
associated to $P$ (see \cite[IV.4.4.2]{Gi}). Then the map defined by  
 $$ \Delta_{Y/X}: \Pic(G) \rightarrow \Br(X), \ \ P \mapsto  \partial_P([Y]) $$ appears in the following short exact sequence (see \cite[Theorem 2.8]{BD})
 \begin{equation} \label{pic-br} \Pic(G) \xrightarrow{\Delta_{X/Y}} \Br(X) \xrightarrow{f^*} \Br(Y) \, . \end{equation}
  
For any $v\in \Omega_k$, the exact sequence \eqref{ses} defines a coboundary map 
$$\partial_P^{k_v}: \ \ \ H^1(k_v, G)\rightarrow H^2(k_v, \Bbb G_m)=\Br(k_v) \, .$$
One can therefore define a pairing 
$$ \delta_v: \ H^1(k_v, G) \times \Pic (G) \to \Br(k_v) \subseteq \Bbb Q/\Bbb Z, \ \ (\sigma_v, P)\mapsto \partial_P^{k_v}(\sigma_v) $$
such that the following diagram
\begin{equation} \label{diag pairings}
 \xymatrix{
 X(k_v)   \ar[d]_{[Y]} & \times & \Br(X)   &  \ar[r]^(.35){ev}  &  \Br(k_v)  \ar[d]^{\id}   \\
 H^1(k_v, G)  & \times &  \Pic(G) \ar[u]_{\Delta_{X/Y}} &  \ar[r]^(.4){\delta_v}  & \Br(k_v) 
 }
\end{equation} 
commutes (see Proposition 2.9 in \cite{CTX}). These pairings induce a pairing 
$$ (\delta_v)_{v\in \Omega_k} : \ \ \ \bigoplus_{v\in \Omega_k} H^1(k_v, G) \times \Pic(G) \rightarrow \Bbb Q/\Bbb Z, \ \ ((\sigma_v)_{v\in \Omega_k}, P)\mapsto \sum_{v\in \Omega_k} \delta_v(\sigma_v, P)\in \Bbb Q/\Bbb Z $$ and a natural exact sequence of pointed sets 
$$ H^1(k, G) \to \bigoplus_{v\in \Omega_k} H^1(k_v, G) \to {\rm Hom} (\Pic(G), \Bbb Q/\Bbb Z) $$
by \cite[Theorem 3.1]{CTX}. Therefore (\ref{equiv}) is equivalent to the fact that $ ([Y](x_v))\in \bigoplus_{v\in \Omega_k} H^1(k_v, G)$ is orthogonal to $\Pic(G)$ for the pairing $(\delta_v)_{v\in \Omega_k}$. The commutative diagram \eqref{diag pairings}, together with \eqref{pic-br}, gives  
$$X({\bf A}_k)^{\ker(f^*)}=\bigcup_{\sigma\in H^1(k, G)} f_\sigma (Y^\sigma ({\bf A}_k)) . $$   Since $\ker(f^*)\subseteq A$, one has 
$$ X({\bf A}_k)^{A}\subseteq X({\bf A}_k)^{\ker(f^*)}=\bigcup_{\sigma\in H^1(k, G)} f_\sigma (Y^\sigma ({\bf A}_k)).  $$
Then the functoriality of the Brauer-Manin pairing implies that 
$$ X({\bf A}_k)^{A}\subseteq \bigcup_{\sigma\in H^1(k, G)} f_\sigma \left(Y^\sigma ({\bf A}_k)^{ f_\sigma^*(A)}\right).  $$

(2) When $G$ is a group of multiplicative type, one obtains that (\ref{equiv}) is equivalent to
$$ \sum_{v\in \Omega_k} \inv_v (\chi \cup [Y]) (x_v) = 0 $$ for all $\chi \in H^1(k, \hat{G})$ by \cite[Theorem 6.3]{D0}.  
Let 
$$ \mathcal{K}_f=\langle \{ \chi\cup [Y] : \  \chi \in H^1(k, \hat{G}) \} \rangle $$
 be the subgroup of $\Br(X)$ generated by elements $\chi\cup [Y]$, where $\cup$ is the cup product 
$$\cup:  H^1(k, \hat{G}) \times H^1(X, G) \rightarrow H^2(X, \Bbb G_m)=\Br(X) .$$ 
Then $$X({\bf A}_k)^{\mathcal{K}_f}=\bigcup_{\sigma\in H^1(k, G)} f_\sigma (Y^\sigma ({\bf A}_k))$$ by \cite[Proposition 3.1]{HaSk}. 
Functoriality of the cup product proves that the following diagram 
$$ \begin{CD}
H^1(k, \hat{G}) \times H^1(X, G) @>{\cup}>> H^2(X, \Bbb G_m)=\Br(X) \\
@V{\id\times f^*}VV @VV{f^*}V \\
H^1(k, \hat{G}) \times H^1(Y, G) @>{\cup}>>  H^2(Y, \Bbb G_m)=\Br(Y)
\end{CD}$$
is commutative. Since $Y\xrightarrow{f} X$ becomes a trivial torsor over $Y$, the above diagram gives $\mathcal{K}_f \subseteq \ker(f^*) $. 
Since $\mathcal{K}_f \subseteq \ker(f^*) \subseteq A$, one has 
$$ X({\bf A}_k)^{A}\subseteq X({\bf A}_k)^{\mathcal{K}_f}=\bigcup_{\sigma\in H^1(k, G)} f_\sigma (Y^\sigma ({\bf A}_k)).  $$
Then the functoriality of the Brauer-Manin pairing implies that 
$$ X({\bf A}_k)^{A}\subseteq \bigcup_{\sigma\in H^1(k, G)} f_\sigma \left(Y^\sigma ({\bf A}_k)^{ f_\sigma^*(A)}\right).  $$
 \end{proof}

\section{Refinement in the toric case}\label{rtc}

In this section, we will refine Theorem \ref{main-general} for torsors under tori.

\begin{thm}\label{tor} Let $f: Y\rightarrow X$ be a torsor under a torus $G$ over a number field $k$. Assume that $X$ is smooth and geometrically integral. Let $ \ker(f^*) \subseteq A \subseteq \Br(X)$ be a subgroup, and for all $\sigma \in H^1(k,G)$, let $B_\sigma \subseteq \Br_1(Y^\sigma) $ be a subgroup such that $$  {f^*}^{-1}\left(\sum_{\sigma\in H^1(k,G)} \psi_\sigma (\widetilde{B_\sigma})\right) \subseteq A  \, ,$$ where $ \Br_a(Y^\sigma)\xrightarrow{\psi_\sigma} \Br_a(Y)$ is the morphism of Definition \ref{brauer-twist} and $\widetilde{B_\sigma}$ is the image of $B_\sigma$ in $\Br_a(Y^\sigma)$.

Then one has
$$ X({\bf A}_k)^A = \bigcup_{\sigma\in H^1(k, G)} f_\sigma \left(Y^\sigma ({\bf A}_k)^{ B_\sigma + f_\sigma^*(A)}\right)  $$ where $Y^\sigma \xrightarrow{f_\sigma}X$ is the twist of $Y\xrightarrow{f} X$ by $\sigma$.
\end{thm}

\begin{proof} Since
$$\bigcup_{\sigma\in H^1(k, G)} f_\sigma \left(Y^\sigma ({\bf A}_k)^{ B_\sigma + f_\sigma^*(A)}\right) \subseteq \bigcup_{\sigma\in H^1(k, G)} f_\sigma \left(Y^\sigma ({\bf A}_k)^{f_\sigma^*(A)}\right) \subseteq  X({\bf A}_k)^A $$ by the functoriality of Brauer-Manin pairing, one only needs to prove the converse inclusion.

Step 1. We first prove the result when $\hat{G}$ is a permutation Galois module. In this case, Shapiro Lemma and Hilbert 90 gives $H^1(K, G)=\{1\}$ for any field extension $K/k$. This implies that
$$ X({\bf A}_k)^A= f \left(Y({\bf A}_k)^{f^*(A)}\right)$$ by the functoriality of Brauer-Manin pairing. 

Let $(x_v)\in X({\bf A}_k)^A$. Then there is $(y_v)\in Y({\bf A}_k)^{f^*(A)}$ such that $(x_v)=f((y_v))$. 

By Proposition 6.10 (6.10.3) in \cite{Sansuc}, the natural sequence
$$ \Br_1(X) \xrightarrow{f^*} \Br_1(Y) \xrightarrow{\lambda} \Br_a(G) $$ 
is exact, and it induces the exact sequence
$$ (f^*)^{-1}(B) \xrightarrow{f^*} B \xrightarrow{\lambda} \Br_a(G) $$ for any subgroup $B\subseteq \Br_1(Y)$. Therefore the following sequence
$$ \Br_a(G)^D \xrightarrow{\lambda^D} B^D \xrightarrow{(f^*)^D} ((f^*)^{-1}(B))^D $$ is exact. Assuming $(f^*)^{-1}(B) \subseteq A$, one has $(f^*)^D((y_v))=0$, where we (abusively) identify $(y_v)$ with its image in $B^D$ via the Brauer-Manin pairing. By the aforementioned exactness, there is $\xi \in \Br_a(G)^D$ such that $\lambda^D(\xi)=(y_v)$. Since $\Sha^1(k,G)=\{1\}$,  Theorem 2 in \cite{Ha08} implies that every element in $\Br_a(G)^D$ is given by an element in $G({\bf A}_k)$ via the Brauer-Manin pairing. Namely, there is $(g_v)\in G({\bf A}_k)$ such that
$$b(y_v)=\lambda(b) (g_v)$$ for all $b\in B$. Then $(g_v)^{-1}\cdot (y_v) \in Y({\bf A}_k)^{B+ f^*(A)}$ by Proposition \ref{torsorbraueraction}, and $(x_v)=f((g_v)^{-1}\cdot (y_v))$.

Step 2. We now prove the case of an arbitrary torus $G$. By Proposition-Definition 3.1 in \cite{CT08}, there is a short exact sequence of tori
$$ 1\rightarrow G \rightarrow T_0 \xrightarrow{q} T_1 \rightarrow 1 \, ,$$
such that $\hat{T_0}$ is a permutation Galois module and $\hat{T_1}$ is a coflasque Galois module. Since
$$ H^3(k, \hat{T_1}) \cong \prod_{v\in \infty_k} H^3(k_v, \hat{T_1}) \cong \prod_{v\in \infty_k} H^1(k_v, \hat{T_1}) =\{1\}$$ 
(see for instance Proposition 5.9 in \cite{HS05}), the map $\Br_1(T_0) \rightarrow \Br_1(G)$ is surjective.

Let $Z\xrightarrow{\rho} X$ be the torsor under $T_0$ defined by $Z := T_0\times^G Y$. We have a morphism of torsors under $G$:
 $$ \begin{CD}
 Y @>{e_0 \times \id_Y}>> T_0\times_k Y @>{\chi}>> Z=T_0 \times^G Y \\
@. @V{p_0}VV  @VV{\theta}V \\
@. T_0 @>>{q}> T_1 \\
\end{CD} $$
where $e_0 \in T_0(k)$ is the unit element, $p_0$ is the projection map and $\theta$ is given as in (\ref{d}). For simplicity, denote by $i := \chi \circ (e_0 \times \id_Y) : Y \to Z$ the composite morphism defined in the previous diagram.

Then Proposition 6.10 (6.10.3) in \cite{Sansuc} gives the following commutative diagram of exact sequences:
$$ \begin{CD}
 \Br_1(T_1) @>{q^*}>> \Br_1(T_0)@>>> \Br_a(G) \\
 @V{\theta^*}VV  @VV{p_0^*}V  @VV{\id}V \\
\Br_1(Z) @>>{\chi^*}> \Br_1(T_0\times_k Y) @>>> \Br_a(G) \, . \\
\end{CD} $$
Since the following sequence
$$ \Br_1(T_0) \xrightarrow{p_0^*} \Br_1(T_0\times_k Y) \xrightarrow{(e_0\times \id_{Y})^*} \Br_a(Y) \rightarrow 1 $$ is exact by Lemma 6.6 in \cite{Sansuc}, the surjectivity of the map $\Br_1(T_0)\rightarrow \Br_1(G)$ implies that the morphism
$$i^* : \ \ \Br_1(Z)\rightarrow \Br_1(Y) $$ is surjective, by a simple diagram chase.

Lemma \ref{q} implies that for any $t\in T_1 (k)$, the composite morphism $\theta^{-1}(t) \to Z\xrightarrow{\rho} X$ is canonically isomorphic to the twist $f_t: Y^{q^{-1}(t)} \rightarrow X$ of $f: Y\rightarrow X$ by the $\Spec(k)$-torsor $q^{-1}(t)$ under $G$.

Denote by $i_t: \theta^{-1}(t) \rightarrow Z$ the closed immersion. Then $f_t=\rho\circ i_t$ for any $t\in T_1(k)$. 

Let $\chi_t$ be the restriction of $\chi$ to $q^{-1}(t) \times_k Y$ for any $t\in T_1(k)$. Then the following diagram
$$ \begin{CD}
 @. q^{-1}(t) \times_k Y @>{\chi_t}>> Y^{q^{-1}(t)}  \\
@. @VV{j_t\times \id_Y}V  @VV{i_t}V \\
Y @>{e_0 \times \id_Y}>> T_0\times_k Y @>{\chi}>> Z \\
@. @VV{p_0}V  @VV{\theta}V \\
G @>>> T_0 @ >{q}>> T_1
\end{CD} $$ is commutative, where $j_t: q^{-1}(t)\rightarrow T_0$ is the closed immersion of the fiber of $q$ at $t$.  
Therefore 
%
 Definition \ref{brauer-twist} implies that we have a commutative triangle:
\[
\xymatrix{
\Br_a(Z) \ar[r]^(.4){i_t^*} \ar[rd]_{i^*} & \Br_a(Y^{q^{-1}(t)}) \ar[d]_{\sim}^{\psi_{q^{-1}(t)}} \\
& \Br_a(Y) \, ,
}
\]
i.e. that 
$\psi_{q^{-1}(t)} \circ i_{t}^* = i^*$.

Let $$B= {i^*}^{-1}\left(\sum_{t\in T_1(k)}  \psi_{q^{-1}(t)} \left(\widetilde{ B_{q^{-1}(t)}}\right)\right) \subset \Br_a(Y)$$
where $\widetilde{ B_{q^{-1}(t)}}$ is the image of $B_{q^{-1}(t)}$ in $\Br_a(Y^{q^{-1}(t)})$ and $\psi_{q^{-1}(t)}$ is given by Definition \ref{brauer-twist} for all $t\in T_1(k)$. 

Since $i^*\circ \rho^* = f^*$, we have 
$${\rho^*}^{-1} (B) ={f^*}^{-1}\left(\sum_{t\in T_1(k)}  \psi_{q^{-1}(t)} \left(\widetilde{ B_{q^{-1}(t)}}\right)\right) \subseteq A \, ,$$ 
hence step 1 applied to the torsor $Z\xrightarrow{\rho} X$ under $T_0$ implies that
\begin{equation} \label{desc qt}
X({\bf A}_k)^A= \rho\left(Z({\bf A}_k)^{B+\rho^*(A)}\right) \, .
\end{equation}

Let $(x_v)\in X({\bf A}_k)^A$. By \eqref{desc qt}, there is $(z_v)\in Z({\bf A}_k)^{B+ \rho^*(A)}$ such that $(x_v)=\rho((z_v))$. 
Since $$i^*\circ \theta^*(\Br_1(T_1))=(e_0\times \id_{Y})^*\circ p_0^* \circ q^* (\Br_1(T_1))=\Br_0(Y) $$ and $i^*(\Br_0(Z))=\Br_0(Y)$, one gets $\theta^*(\Br_1(T_1))\subseteq \Br_0(Z) + B$ (by construction, $B$ contains $\ker(i^* : \Br_1(Z) \to \Br_1(Y))$). Functoriality of the Brauer-Manin pairing now gives 
$$\theta ((z_v))\in T_1({\bf A}_k)^{\Br_1(T_1)}\, .$$
 By Proposition \ref{connected-groups}, there are $\alpha\in T_1(k)$ and $(\beta_v)\in T_0({\bf A}_k)^{\Br_1(T_0)}$ such that  $\theta ((z_v))= \alpha \cdot q(\beta_v)$. 

Therefore $ (\beta_v)^{-1} \cdot (z_v) \in \theta^{-1}(\alpha)$, hence $(\beta_v)^{-1} \cdot (z_v) \in Z({\bf A}_k)^{B+ \rho^*(A)}$. 

Since $i^* : \Br_1(Z) \to \Br_1(Y)$ is surjective, one has 
$$\psi_{q^{-1}(\alpha)} \circ i_\alpha^*(\widetilde{B}) = i^* (\widetilde{B}) = \sum_{t\in T_1(k)}  \psi_{q^{-1}(t)} \left(\widetilde{ B_{q^{-1}}}\right)  \supseteq \psi_{q^{-1}(\alpha)}\left(\widetilde{B_{q^{-1}(\alpha)}}\right) \, ,$$ where $\widetilde{B}$ is the image of $B$ in $\Br_a(Z)$. It implies that $i_\alpha^*(B)+ \Br_0(\theta^{-1}(\alpha)) \supseteq B_{q^{-1}(\alpha)}$ by Lemma \ref{twist-isom}, and 
$$ (\beta_v)^{-1} \cdot (z_v) \in \left[\theta^{-1}(\alpha) ({\bf A}_k)\right]^{i_\alpha^*(B)+ (i_\alpha^*\circ \rho^*)(A)} \subseteq \left[\theta^{-1}(\alpha) ({\bf A}_k)\right]^{B_{q^{-1}(\alpha)}+ (i_\alpha^*\circ \rho^*)(A)} $$ as desired. 
\end{proof}

The first part of the following result is also proved in Theorem 1.7 of \cite{Wei}.

\begin{cor} \label{algebraic} Let $X$ be a smooth and geometrically integral variety. If $f: Y\rightarrow X$ is a torsor under a torus $G$ over a number field $k$, then
$$ X({\bf A}_k)^{\Br_1(X)} = \bigcup_{\sigma\in H^1(k, G)} f_\sigma \left(Y^\sigma ({\bf A}_k)^{\Br_1(Y^\sigma)}\right) $$
and 
$$ X({\bf A}_k)^{\Br} = \bigcup_{\sigma\in H^1(k, G)} f_\sigma \left(Y^\sigma ({\bf A}_k)^{\Br_1(Y^\sigma)+f_\sigma^*(\Br(X))}\right)  . $$
\end{cor}

\begin{proof} To get the first equality, apply Theorem \ref{tor} to $A=\Br_1(X)$ and $B_\sigma=\Br_1(Y^\sigma)$ for each $\sigma \in H^1(k, G)$. Since $\Pic (G_{\bar k})=0$, Proposition 6.10 in \cite{Sansuc} gives
$${f^*}^{-1}\left(\sum_{\sigma\in H^1(k,G)} \psi_\sigma \left(\widetilde{B_\sigma}\right)\right) \subseteq {f^*}^{-1}(\Br_a(Y)) \subseteq \Br_1(X)=A \, ,$$ 
as required.

The second equality follows from Theorem \ref{tor} by taking $A=\Br(X)$ and $B_\sigma=\Br_1(Y^\sigma)$ for each $\sigma \in H^1(k, G)$. 
\end{proof}

\section{An application} \label{aa}
In this section, we apply the previous results to study the necessary conditions for a connected linear algebraic group to satisfy strong approximation with Brauer-Manin obstruction.  

When $X$ is affine, the set $X(k)$ is discrete in $X({\bf A}_k)$ by the product formula. Therefore if such an $X$ satisfies strong approximation off $S$, then $\prod_{v\in S} X(k_v)$ is not compact. However this necessary condition for strong approximation is no longer true for strong approximation with Brauer-Manin obstruction if $\Br(X)/\Br(k)$ is not finite. For example, a torus $X$ always satisfies strong approximation with Brauer-Manin obstruction off $\infty_k$, $X$ being anisotropic over $k_\infty$ or not: see \cite[Theorem 2]{Ha08}. When $X$ is a semi-simple linear algebraic group, the necessary and sufficient condition for $X$ to satisfy strong approximation with Brauer-Manin obstruction is given by Proposition 6.1 in \cite{CX2}. In this section, we extend this result to a general connected linear algebraic group.    

The following lemma explains that strong approximation with Brauer-Manin obstruction for a general connected linear algebraic group can be reduced to the reductive case. 
\begin{lem} \label{unip-equ} Let $G$ be a connected linear algebraic group over a number field $k$.

If $\pi: G\to G^{\rm{red}}$ is the quotient map, then $G^{\rm{red}}({\bf A}_k)^{\Br_1(G^{\rm{red}})} = \pi \left(G({\bf A}_k)^{\Br_1(G)} \right)$.

In particular, for any finite subset $S$ of $\Omega_k$, $G$ satisfies strong approximation with respect to $\Br_1(G)$ off $S$ if and only if $G^{\rm{red}}$ satisfies strong approximation with respect to $\Br_1(G^{\rm{red}})$ off $S$.
\end{lem}

\begin{proof} By applying Lemma \ref{br} for $k$ and $\bar{k}$, one obtains that $ \pi^*(\Br_1(G^{\rm{red}}))=\Br_1(G)$. 
The first part follows from Theorem \ref{main-general} and Proposition 6 of \S 2.1 of  Chapter III in \cite{Ser}. 

Suppose $G$ satisfies strong approximation with respect to $\Br_1(G)$ off $S$. For any open subset $$M=\prod_{v\in S} G^{\rm{red}}(k_v) \times \prod_{v\not\in S} M_v $$ of $G^{\rm{red}}({\bf A}_k)$ such that $M \cap \left[G^{\rm{red}}({\bf A}_k)^{\Br_1(G^{\rm{red}})}\right] \neq \emptyset$, one has that 
$$ \pi^{-1}(M) = \prod_{v\in S} G(k_v) \times \prod_{v\not\in S} \pi^{-1}(M_v) $$ with 
$\pi^{-1}(M) \cap G({\bf A}_k)^{\Br_1(G)} \neq \emptyset$ by the first part. Then by assumption there is $x\in G(k) \cap \pi^{-1}(M)$. It implies that $\pi(x)\in M\cap G^{\rm{red}}(k)$, as required. 

Conversely, suppose $G^{\rm{red}}$ satisfies strong approximation with respect to $\Br_1(G^{\rm{red}})$ off $S$. For any open subset $$N=\prod_{v\in S} G(k_v) \times \prod_{v\not\in S} N_v $$  of $G({\bf A}_k)$ such that $N\cap G({\bf A}_k)^{\Br_1(G)} \neq \emptyset$, we have 
$$ \pi (N) = \prod_{v\in S} G^{\rm{red}}(k_v) \times \prod_{v\not\in S} \pi (N_v) $$ 
and this set is an open subset of $G^{\rm{red}}({\bf A}_k)$, with $\pi(M) \cap \left[G^{\rm{red}}({\bf A}_k)^{\Br_1(G^{\rm{red}})}\right] \neq \emptyset$: here we use Proposition 6 of \S 2.1 of  Chapter III in \cite{Ser}, Proposition 6.5 in Chapter 6 of \cite{PR} and the functoriality of Brauer-Manin pairing. Then by assumption there is $y\in G^{\rm{red}}(k) \cap \pi(N)$. Using Proposition 6 of \S 2.1 of  Chapter III in \cite{Ser} one more time, one concludes that $\pi^{-1}(y)$ is isomorphic to $R_u(G)$ as an algebraic variety, hence it satisfies strong approximation off $S$. Since 
$$ \pi^{-1}(y)\cap N= \prod_{v\in S} \pi^{-1}(y)(k_v) \times \prod_{v\not\in S} (\pi^{-1}(y)(k_v) \cap N) \neq \emptyset , $$
there is $z\in \pi^{-1}(y)(k) \cap N \subset G(k)\cap N$, as desired. 
\end{proof}

The main result of this section is the following statement:

\begin{thm} \label{red-semi} Let $G$ be a connected linear algebraic group over a number field $k$ and let $G^{\rm{qs}}:=G/R(G)$, where $R(G)$ is the solvable radical of $G$. If $\pi: G\to G^{\rm{qs}}$ is the quotient map, then $$G^{\rm{qs}}({\bf A}_k)^{\Br_1(G^{\rm{qs}})} = \pi \left(G({\bf A}_k)^{\Br_1(G)} \right) \cdot G^{\rm{qs}}(k) \, . $$ 

 In particular, if $G$ satisfies strong approximation with respect to $\Br_1(G)$ off a finite subset $S$ of $\Omega_k$, then $G^{\rm{qs}}$ satisfies strong approximation with respect to $\Br_1(G^{\rm{qs}})$ off $S$.
\end{thm}
\begin{proof} For the first part, by functoriality of the Brauer-Manin pairing, one only needs to prove that 
$$G^{\rm{qs}}({\bf A}_k)^{\Br_1(G^{\rm{qs}})} \subseteq \pi \left(G({\bf A}_k)^{\Br_1(G)} \right) \cdot  G^{\rm{qs}}(k) \, .$$ 

By Lemma \ref{unip-equ}, we can assume that $G$ is reductive. Then $R(G)$ is a torus contained in the center of $G$ (see Theorem 2.4 in Chapter 2 of \cite{PR}) and $\pi: G\to G^{\rm{qs}}$ is a torsor under $R(G)$. By Corollary \ref{algebraic}, for any $(x_v)\in G^{\rm{qs}}({\bf A}_k)^{\Br_1(G^{\rm{qs}})}$, there are $\sigma\in H^1(k, R(G))$ and $(y_v) \in G^{\sigma}({\bf A}_k)^{\Br_1(G^\sigma)}$ such that $(x_v)=\pi_\sigma ((y_v))$.  Since $G^\sigma(k)\neq \emptyset $ by Corollary 8.7 in \cite{Sansuc} (see also Theorem 5.2.1 in \cite{Sko}),  there is $\gamma\in G^{\rm{qs}}(k)$ such that $\partial(\gamma)=\sigma$, where $\partial$ is the coboundary map in the following exact sequence in Galois cohomology:
$$ 1\to R(G)(k) \to G(k) \to G^{\rm{qs}}(k) \xrightarrow{\partial} H^1(k, R(G)) \to H^1(k, G) \, .$$

In addition, the choice of an element $\bar{\gamma} \in G(\bar{k})$ such that $\pi(\bar{\gamma}) = \gamma$ defines a commutative diagram defined over $k$:
$$ \begin{CD}
 G^\sigma  @>{\bar{\gamma} \, \cdot}>{\sim}>  G \\
@V{\pi_\sigma}VV  @VV{\pi}V \\
G^{\rm{qs}}  @>{\gamma \, \cdot}>{\sim}>  G^{\rm{qs}} \\
\end{CD} $$
(see for instance Example 2 of p.20 in \cite{Sko}).  This implies that $$ \pi_\sigma \left(G^\sigma({\bf A}_k)^{\Br_1(G^\sigma)} \right) = \pi \left(G({\bf A}_k)^{\Br_1(G)} \right) \cdot \gamma \, ,$$ as desired.

Suppose now that $G$ satisfies strong approximation with respect to $\Br_1(G)$ off $S$. For any open subset $$M=\prod_{v\in S} G^{\rm{qs}}(k_v) \times \prod_{v\not\in S} M_v $$ of $G^{\rm{qs}}({\bf A}_k)$ such that $M \cap G^{\rm{qs}}({\bf A}_k)^{\Br_1(G^{\rm{qs}})} \neq \emptyset$, the first part implies that there is $g \in G^{\rm{qs}}(k)$ such that 
$$ \pi^{-1}(M\cdot g) = \prod_{v\in S} G(k_v) \times \prod_{v\not\in S} \pi^{-1}(M_v\cdot g)\, ,$$ 
with $\pi^{-1}(M\cdot g) \cap G({\bf A}_k)^{\Br_1(G)} \neq \emptyset$. Since $G$ satisfies strong approximation with algebraic Brauer-Manin obstruction off $S$, there exists $x\in G(k) \cap \pi^{-1}(M \cdot g)$. This implies that $\pi(x)\cdot g^{-1}\in M\cap G^{\rm{qs}}(k)$ as required.  \end{proof}

\begin{cor} \label{iff} Let $G$ be a connected linear algebraic group over a number field $k$ and let $S$ a finite subset of $\Omega_k$ containing $\infty_k$. Then $G$ satisfies strong approximation with respect to $\Br_1(G)$ off $S$ if and only if $\prod_{v\in S} G'(k_v)$ is not compact for any non-trivial simple factor $G'$ of the semi-simple part $G^{ss}$ of $G$. 
\end{cor}
\begin{proof} By Theorem 2.3 and Theorem 2.4 of Chapter 2 in \cite{PR}, the quotient map $$G^{\rm{red}}\to G/R(G)=G^{\rm{qs}}$$ induces an isogeny $G^{ss} \to G^{\rm{qs}}$. 
One side follows from Corollary 3.20 in \cite{D}. The other side follows from Theorem \ref{red-semi} and Proposition 6.1 in \cite{CX2}.
\end{proof}

\begin{rem} All the results in this section involve the group $\Br_1(G)$, and they remain true with $\Br_1(G)$ replaced by $\Br(G)$. Indeed, there is a sufficiently large subset $S$ of $\Omega_k$ containing $\infty_k$ such that $\prod_{v\in S} G'(k_v)$ is not compact for any non-trivial simple factor $G'$ of $G^{ss}$, therefore Corollary 3.20 in \cite{D}, Proposition 2.6 in \cite{CTX} and the functoriality of Brauer-Manin pairing gives the following inclusions:
$$ G({\bf A}_k)^{\Br_1(G)} = \overline{G(k) \cdot \rho (\prod_{v\in S} G^{scu}(k_v))} \subseteq  G({\bf A}_k)^{\Br(G)} \subseteq  G({\bf A}_k)^{\Br_1(G)} \, , $$
where $G^{scu}=G^{sc}\times_{G^{\rm{red}}} G$ with the projection map $G^{scu}\xrightarrow{\rho} G$ and $G^{sc}$ is the simply connected covering of $G^{ss}$. In particular, we have $G({\bf A}_k)^{\Br(G)} = G({\bf A}_k)^{\Br_1(G)}$.
\end{rem}

\section{Comparison I, $X({\bf A}_k)^{\rm{desc}} \subseteq X({\bf A}_k)^{\textup{\'et}, \Br}$ }\label{cI}

Let $Y\xrightarrow{f} X$ be a left torsor under a linear algebraic group $G$ over a number field $k$. The fundamental problem to define the descent obstruction for strong approximation with respect to $Y\xrightarrow{f} X$ is to decide whether the set 
$$ X({\bf A}_k)^f = \left\{(x_v)\in X({\bf A}_k): ([Y](x_v))\in {\rm{Im}} \left(H^1(k, G) \rightarrow \prod_v H^1(k_v, G)\right) \right\} = \bigcup_{\sigma\in H^1(k, G)} f_\sigma (Y^\sigma ({\bf A}_k)) $$ is closed or not in $X({\bf A}_k)$. We already know that this is true when $G$ is either connected or a group of multiplicative type, by Theorem \ref{main-general}. For a general linear algebraic group $G$, this result is proved by Skorobogatov in Corollary 2.7 of \cite{ Sk1}, when $X$ is assumed to be proper over $k$. The proof depends on Proposition 5.3.2 in \cite{Sko} or Proposition 4.4 in \cite{HSk}, which  are not true for open varieties, as explained in the following example.

\begin{exa} The short exact sequence of linear algebraic groups 
$$ 1 \rightarrow \mu_2 \rightarrow \Bbb G_m \xrightarrow{f} \Bbb G_m \rightarrow 1 \, ,$$ 
where $f(x)=x^2$, can be viewed as torsor over $\G_m$ under $\mu_2$.  For any $\sigma \in H^1(k, \mu_2)\cong k^\times/(k^\times)^2$, the twist $\Bbb G_m^\sigma$ of $\Bbb G_m$ by $\sigma$ is given by the equation $x=a_\sigma y^2$ in $\Bbb G_m\times_k \Bbb G_m$, where $a_\sigma$ is an element in $k^\times$ representing the class $\sigma$ by the above isomorphism. It is clear that $\Bbb G_m^\sigma \cong \Bbb G_m$ as varieties over $k$, hence it always contains adelic points. 
\end{exa}

We use the same definition of an integral model as in \cite{LX}.

\begin{definition} Let $X$ be a variety over a number field $k$ and let $S$ be a finite subset of $\Omega_k$ containing $\infty_k$. An integral model of $X$ over $O_S$ is a faithfully flat separated $O_S$-scheme $\mathcal{X}_S$ of finite type such that $\mathcal{X}_S\times_{O_S} k \cong X$. 
\end{definition} 

The replacement for Proposition 5.3.2 in \cite{Sko} or Proposition 4.4 in \cite{HSk} is the following proposition:

\begin{prop}\label{finite} Let $X$ be a variety over a number field $k$ and let $S$ be a finite subset of $\Omega_k$ containing $\infty_k$. Fix an integral model $\mathcal{X}_S$ of $X$ over $O_S$. If $Y\xrightarrow{f} X$ is a left torsor under a linear algebraic group $G$ over $k$, then the set 
$$ \left\{ [\sigma] \in H^1(k, G): \ f_{\sigma} (Y^\sigma ({\bf A}_k))\cap \left[ \prod_{v\in S}X(k_v)  \times \prod_{v\not\in S} \mathcal{X}_S(O_v) \right]\neq \emptyset \right\} $$ is finite.
\end{prop}
\begin{proof} It follows from the same argument as the proof of Proposition 4.4 in \cite{HSk}.
\end{proof}

 One can now extend Corollary 2.7 in \cite{ Sk1} to open varieties by using the above replacement for Proposition 4.4 in \cite{HSk}. 

\begin{prop}\label{closed} Let $X$ be a (not necessarily proper) variety over a number field $k$. If $Y\xrightarrow{f} X$ is a left torsor under a linear algebraic group $G$ over $k$, then the set
$X({\bf A}_k)^f$ is closed in $X({\bf A}_k)$.
\end{prop}
\begin{proof} Take an integral model $\mathcal{X}_{S_0}$ of $X$ over $O_{S_0}$, where $S_0$ is a finite subset of $\Omega_k$ containing $\infty_k$.  Then 
$$ \left\{ \prod_{v\in S} X(k_v) \times \prod_{v\in \Omega_k\setminus S} \mathcal{X}_{S_0} (O_v) \right\}_{S} $$ is an open covering of $X({\bf A}_k)$ (see Theorem 3.6 in \cite{Conrad}), where $S$ runs through all finite subsets of $\Omega_k$ containing $S_0$. By Proposition \ref{finite} and Corollary 2.5 in \cite{Sk1}, the set 
$$ X({\bf A}_k)^f \cap \left[\prod_{v\in S} X(k_v) \times \prod_{v\in \Omega_k\setminus S} \mathcal{X}_{S_0} (O_v)\right] $$
is closed in $\prod_{v\in S} X(k_v) \times \prod_{v\in \Omega_k\setminus S} \mathcal{X}_{S_0} (O_v)$, therefore the set $X({\bf A}_k)^f$ is closed in $X({\bf A}_k)$.
\end{proof}

Applying Proposition \ref{finite}, one can also extend Lemma 2.2 and Theorem 1.1 in \cite{Sk1} to open varieties.
For any variety over a number field $k$, and following \cite{Sk1}, we write
$$ X({\bf A}_k)^{\text{desc}} =\bigcap_{Y\xrightarrow{f} X} X({\bf A}_k)^f \, ,$$
where $Y\xrightarrow{f} X$ runs through all torsors under all linear algebraic groups over $k$ (see also \S \ref{intro}.).

\begin{lem} \label{rep2.2} Let $X$ be a (not necessarily proper) variety and let $Y\rightarrow X$ be a torsor over a number field $k$. For any $(P_v)\in  X({\bf A}_k)^{\rm{desc}}$, there is a twist $Y'\rightarrow X$ of $Y\rightarrow X$ such that the following property holds:

For any surjective $X$-torsor morphism $Z\rightarrow Y'$ (see Definition 2.1 in \cite{Sk1}), there is a twist $Z'\rightarrow Y'$ of $Z\rightarrow Y'$ such that $(P_v)$ lies in the image of $Z'({\bf A}_k)$. \end{lem}

\begin{proof} There are a finite subset $S_0$ of $\Omega_k$ containing $\infty_k$ and an integral model $\mathcal{X}_{S_0}$ over $O_{S_0}$ such that 
$$ (P_v) \in \prod_{v\in S_0} X(k_v) \times \prod_{v\in \Omega_k\setminus S_0} \mathcal{X}_{S_0} (O_v) $$ (see for instance Theorem 3.6 in \cite{Conrad}), hence Proposition \ref{finite} implies that there are only finitely many twists of a given torsor over $X$ such that $(P_v)$ lifts as an adelic point of this torsor. As pointed out in the proof of Lemma 2.2 in \cite{Sk1}, the finite combinatorics in the first part of the proof of Proposition 5.17 in \cite{St} are still valid.  It concludes the proof.
\end{proof}

\begin{prop} \label{dd}  Let $X$ be a (not necessarily proper) variety over a number field $k$. If $Y\xrightarrow{f} X$ is a left torsor under a finite group scheme $F$ over $k$, then
$$ X({\bf A}_k)^{\rm{desc}} =\bigcup_{\sigma\in H^1(k, F)} f_{\sigma} \left(Y^{\sigma}({\bf A}_k)^{\rm{desc}}\right)  . $$
\end{prop}

\begin{proof}  One only needs to modify the proof of Theorem 1.1 in \cite{Sk1} by replacing Lemma 2.2 in \cite{Sk1} with Lemma \ref{rep2.2}, Corollary 2.7 in \cite{Sk1} with Proposition \ref{closed}. Moreover, since $f$ is finite, the induced map $Y({\bf A}_k) \xrightarrow{f} X({\bf A}_k)$ is topologically proper by Proposition 4.4 in \cite{Conrad}. This implies that $f^{-1}((P_v))$ is compact. \end{proof}

Recall that, following \cite{P}, one can define for any variety $X$ over a number field $k$, the set 
$$ X({\bf A}_k)^{\textup{\'et}, \Br}= \bigcap_{Y\xrightarrow{f} X} \bigcup_{\sigma\in H^1(k, F)} f_\sigma (Y^\sigma({\bf A}_k)^{\Br}) \, ,$$ 
where $Y\xrightarrow{f} X$ runs over all torsors under all finite groups $F$ over $k$ (see \S \ref{intro}).  Since the induced map $Y({\bf A}_k) \xrightarrow{f} X({\bf A}_k)$ is topologically closed for any finite morphism $Y\xrightarrow{f} X$ by Proposition 4.4 in \cite{Conrad}, one concludes that $X({\bf A}_k)^{\textup{\'et}, \Br}$ is closed in $ X({\bf A}_k)$ by the same argument as in Proposition \ref{closed}. 

\begin{cor} \label{oneside} If $X$ is a smooth quasi-projective variety over a number field $k$, then
$$   X({\bf A}_k)^{\rm{desc}} \subseteq X({\bf A}_k)^{\textup{\'et}, \Br} \subseteq  X({\bf A}_k)^{\Br} . $$
\end{cor}

\begin{proof} One only needs to show that  $ X({\bf A}_k)^{\text{desc}} \subseteq X({\bf A}_k)^{\textup{\'et}, \Br}$. For any torsor $Y\xrightarrow{f} X$ under a finite group scheme $F$,  Proposition \ref{dd} gives the equality
$$ X({\bf A}_k)^{\text{desc}} =\bigcup_{\sigma\in H^1(k, F)} f_{\sigma} \left(Y^{\sigma}({\bf A}_k)^{\text{desc}}\right) \, .$$
Since $X$ is quasi-projective, $Y^\sigma$ is quasi-projective as well. By a theorem of Gabber (see \cite{de}), one has 
$$Y^{\sigma}({\bf A}_k)^{\text{desc}}\subseteq Y^{\sigma}({\bf A}_k)^{\Br} $$ (see the proof of Lemma 2.8 in \cite{Sk1}) and the result follows.
\end{proof}

\section{Comparison II, $X({\bf A}_k)^{\textup{\'et}, \Br}\subseteq X({\bf A}_k)^{\rm{desc}}$} \label{CII}

In this section, we prove the inclusion $X({\bf A}_k)^{\textup{\'et}, \Br}\subseteq X({\bf A}_k)^{\rm{desc}}$ for open varieties, which implies Theorem \ref{inteq}. The strategy of proof is the same as in \cite{D09}.

The second named author would like to thank Laurent Moret-Bailly warmly for finding a mistake and for suggesting the following alternative proof of Lemma 4 in \cite{D09} (which already appeared in \cite{Dth}). The statement of this lemma is correct, but the proof in \cite{D09} uses a result of Stoll (see \cite{St}) that is not. Note that in contrast with \cite{D09}, all torsors (unless explicitely mentioned) are assumed to be left torsors.

\begin{lem} \label{Stoll}
Let $X$ be a smooth geometrically connected $k$-variety. Let $(P_v) \in X({\bf A}_k)^{\et, \Br}$ and let $Z \xrightarrow{g} X$ be a torsor under a finite $k$-group $F$.

Then there are a cocycle $\sigma \in Z^1(k, F)$ and a connected component $X'$ of $Z^\sigma$ over $k$ such that the restriction of $g_\sigma$ to $X'$ is a torsor $X' \rightarrow X$ under the stabilizer $F'$ of $X'$ for the action of $F^\sigma$, and the point $(P_v)$ lifts to a point $(Q'_v) \in X'({\bf A}_k)^{\textup{Br}}$.

In particular, $X'$ is geometrically integral.
\end{lem}
\begin{proof}
By assumption, the point
$(P_v)$ lifts to some point $(Q_v) \in Z^{\sigma}({\bf A}_k)^{\textup{Br}}$ for some cocycle $\sigma$ with values in $F$. Since $Z^{\sigma}$ is smooth, $Z^\sigma$ is a disjoint union of connected components over $k$. By Proposition 3.3 in \cite{LX}, there is a $k$-connected component $X'$ of $Z^{\sigma}$ such that $(Q_v)_{v\not\in \Xi} \in P_{\Xi} (X'({\bf A}_k)^{\textup{Br}})$, where $\Xi$ is the set of all complex places of $k$, ${\bf A}_k^{\Xi}$ is the ring of adeles without $\Xi$-components and $P_{\Xi}$ is the projection from $X'({\bf A}_k)$ to $X'({\bf A}_k^{\Xi})$. Since for $v\in \Xi$, $Z^\sigma \times_k k_v$ is a trivial torsor under the finite constant group scheme $F^\sigma\times_k k_v$, we have $g_\sigma (X'(k_v))=X(k_v)$ for all $v\in\Xi$. Hence one can assume that $Q_v\in X'(k_v)$ for $v\in \Xi$, so that we have $(Q_v) \in X'({\bf A}_k)^{\textup{Br}}$. 

Since $X'$ is connected and $X'({\bf A}_k) \neq \emptyset$, the proof of Lemma 5.5 in \cite{St} implies that $X'$ is geometrically connected. Eventually, $X'$ being geometrically connected guarantees that the variety $X'$ is an $X$-torsor under the stabilizer $F'$ of $X'$ in $F^{\sigma}$.
\end{proof}

Let us continue the proof of the aforementioned inclusion. Let $X$ be a smooth and geometrically integral $k$-variety, and $(P_v) \in X({\bf A}_k)^{\textup{\'et}, \Br}$. We need to prove that $(P_v) \in X({\bf A}_k)^{\rm{desc}}$.

For a linear algebraic group $G$ over $k$, one has the following short exact sequence of algebraic groups over $k$:
$$ 1\rightarrow H\rightarrow G \rightarrow F \rightarrow 1 \, ,$$
where $H$ is the connected component of $G$ and $F$ is finite over $k$. This induces the following diagram of short exact sequences
\begin{displaymath}
	\xymatrix{
1 \ar[r] & H \ar[r] \ar[d] & G \ar[r] \ar[d] & F \ar[r] \ar[d] & 1 \\
1 \ar[r] & T \ar[r] & G' \ar[r] & F \ar[r] & 1 
}\end{displaymath}
where $T$ denotes the maximal toric quotient of $H$ and $G'$ is the quotient of $G$ by the kernel of $H\rightarrow T$.

Let  $Y \rightarrow X$ be a torsor under $G$ and let $Z \rightarrow X$ be the push-forward of $Y \rightarrow X$ by the morphism $G \rightarrow F$, which is a torsor under $F$. If $\sigma \in Z^1(k, F)$ is a $1$-cocycle given by Lemma \ref{Stoll} applied to the torsor $Z \rightarrow X$ and to the point $(P_v)$, we want to show that the cocycle $\sigma \in Z^1(k, F)$ lifts to a cocycle $\tau \in Z^1(k, G)$, as in Proposition 5 in \cite{D09}. The obstruction to lift $\sigma$ to a cocycle in $Z^1(k,G)$ gives a natural cohomology class $\eta_\sigma \in H^2(k,\kappa_\sigma)$ by  (5.1) in \cite{FSS} (see also (7.7) in \cite{B}), where $\kappa_\sigma$ is a natural $k$-kernel on $H_{\bar{k}}$ associated to $\sigma$. Lemma 6 in \cite{D09} implies that there is a canonical map $H^2(k, \kappa_\sigma) \to H^2(k, T^\sigma)$ such that the class $\eta_\sigma$ is neutral if and only if its image $\eta'_\sigma \in H^2(k, T^\sigma)$ is zero.

We now apply the open descent theory and the extended type developed by Harari and Skorobogatov in \cite{HaSk} to establish the analogue of Lemma 7 in \cite{D09} for open varieties.  
As in the proof of \cite{D09}, the torsor $Y \to Z$ under $H$ induces a torsor $ W \xrightarrow{\varpi} Z$ under $T$ by the natural map $H^1(Z,H)\to H^1(Z,T)$. Instead of using the  type of the torsor $\varpi$ that was used in \cite{D09}, we consider the so-called "extended type" of the torsor $\varpi$ that was introduced by Harari and Skorobogatov (see  Definition 8.2 in \cite{HaSk}). For a variety $Z$ over $k$, let $KD'(Z)$ denote the complex of Galois modules $[\overline{k}(Z)^* / \overline{k}^* \to \Div(Z_{\bar k})]$ in the derived category $D^b_\et(k)$ of bounded complexes of \'etale sheaves over $\Spec(k)$. One can associate to the torsor  $W \xrightarrow{\varpi} Z$ under $T$ a canonical morphism in this derived category
$$\lambda_W : \widehat{T} \to KD'(Z) \, , $$
called the extended type of $\varpi$. This induces a morphism in the derived category of bounded complexes of abelian groups
$$\lambda_W^\sigma : \widehat{T}^\sigma \to KD'(Z^\sigma) $$ for the above $\sigma \in Z^1(k,F)$. 

\begin{lem} \label{twist}
The morphism $\lambda_W^\sigma : \widehat{T}^\sigma \to KD'(Z^\sigma)$ is a morphism in the derived category of bounded complexes of \'etale sheaves over $\Spec(k)$.
\end{lem}

\begin{proof}
The natural left actions of $F$ on both $T$ and $Z$ induces right actions of $F$ on $\widehat{T}$ and on $KD'(Z)$.

We first prove that the morphism $\lambda_W$ is $F$-equivariant for those actions. 

Let $f \in F(\overline{k})$. We denote by $f_{Z} : Z_{\bar k} \to Z_{\bar k}$ the morphism of $\bar{k}$-varieties defined by $z \mapsto f \cdot z$. This morphism induces a natural morphism in the derived category $f_{Z}^* : KD'(Z_{\bar k}) \to KD'(Z_{\bar k})$. Similarly, the element $f$ defines a natural morphism of $\bar{k}$-tori $f_T : T_{\bar k} \to T_{\bar k}$ such that $f_T(t) := g t g^{-1}$, where $g \in G'(\overline{k})$ is any point lifting $f \in F(\overline{k})$. This morphism $f_T$ induces a morphism of abelian groups $\widehat{f_T} : \widehat{T} \to \widehat{T}$ such that $\widehat{f_T}(\chi) := \chi \circ f_T$.

One needs to prove that the following diagram
\begin{displaymath}
\xymatrix{
\widehat{T} \ar[r]^{\lambda_{W_{\bar k}} \ \ \ \ } \ar[d]_{\widehat{f_T}} & KD'(Z_{\bar k}) \ar[d]^{f_{Z}^*} \\
\widehat{T} \ar[r]_{\lambda_{W_{\bar k}} \ \ \ \ } & KD'(Z_{\bar k})
}
\end{displaymath}
is commutative.

Let $f_{T, *} W_{\bar k}$ be the push-forward of the torsor $W_{\bar k} \to Z_{\bar k}$ under $T_{\bar k}$ by the $\bar k$-morphism $ T_{\bar k} \xrightarrow{f_T} T_{\bar k}$ and let $f_{Z}^* W_{\bar{k}}$ be the pullback of the torsor $W_{\bar{k}} \to Z_{\bar{k}}$ under $T_{\bar{k}}$ by the $\bar{k}$-morphism $f_Z : Z_{\bar{k}} \to Z_{\bar{k}}$.
Then functoriality of the extended type gives:
$$f_Z^* \circ \lambda_{W_{\bar{k}}} = \lambda_{f_Z^* W_{\bar{k}}}  \, \, \ \ \textup{and} \ \ \, \, \lambda_{f_{T,*} W_{\bar{k}}} = \lambda_{W_{\bar{k}}} \circ \widehat{f_T} \, .$$

To prove the required commutativity $f_Z^* \circ \lambda_{W_{\bar{k}}} = \lambda_{W_{\bar{k}}} \circ \widehat{f_T}$, it is enough to show that the torsors $f_Z^*  W_{\bar{k}} \to  Z_{\bar{k}}$ and $f_{T, *}  W_{\bar{k}} \to  Z_{\bar{k}}$ under $T_{\bar{k}}$ are isomorphic. Indeed, we have the following commutative diagram
\begin{displaymath}
\xymatrix{
T_{\bar{k}} \times W_{\bar{k}} \ar[d]_{\varpi \circ p_W} \ar[r]^{ \ \ \ g} & W_{\bar{k}} \ar[d]^{\varpi} \\
Z_{\bar{k}}\ar[r]_{f_Z} & Z_{\bar{k}} \, ,
}
\end{displaymath}
where $p_W$ denotes the projection on $W_{\bar{k}}$ and the morphism $g$ is defined by $(t,w) \mapsto (t g) \cdot w$. This diagram induces a natural $Z_{\bar{k}}$-morphism $\phi : T_{\bar{k}} \times W_{\bar{k}} \to f_Z^* W_{\bar{k}}$. Consider now the right action of $T_{\bar{k}}$ on $T_{\bar{k}} \times W_{\bar{k}}$ defined by $(s,w) \cdot t := (s f_T(t),t^{-1} \cdot w) = (s g t g^{-1},t^{-1} \cdot w)$. Then the morphism $\phi$ is $T_{\bar{k}}$-invariant under this action, hence it induces a $Z_{\bar{k}}$-morphism $\psi : f_{T,*} W_{\bar{k}} \to f_Z^* W_{\bar{k}}$. One can check by a simple computation that $\psi$ is $T_{\bar{k}}$-equivariant, i.e. that $\psi$ is a morphism of (left) torsors over $Z_{\bar{k}}$ under $T_{\bar{k}}$. It concludes the proof of the required commutativity, hence the morphism $\lambda_W$ is $F$-equivariant.

By definition of the twists $T^\sigma$ and $Z^\sigma$, the fact that $\lambda_W$ is $F$-equivariant implies that the morphism $\lambda_W^\sigma$ is Galois equivariant, i.e. that $\lambda_W^\sigma$ is a morphism in the derived category of bounded complexes of \'etale sheaves over $\Spec(k)$.
\end{proof}

By Proposition 8.1 in \cite{HaSk}, there is a natural exact sequence of abelian groups
$$H^1(k, T^\sigma) \to H^1(X', T^\sigma) \xrightarrow{\lambda} \Hom_k(\widehat{T^\sigma}, KD'(X')) \xrightarrow{\partial} H^2(k, T^\sigma) $$
where the map $\lambda$ is the extended type. Let $\lambda_\sigma '=\psi^* \circ \lambda_W^\sigma$, where $\psi : X' \to W$ is the inclusion of the $k$-connected component given by Lemma \ref{Stoll}, and $KD'(Z^\sigma)\xrightarrow{\psi^*} KD'(X')$ is the map induced by $\psi$.

The following lemma, which is an analogue of Lemma 8 in \cite{D09}, is a crucial step for proving the main result of this section. We give here a more conceptual proof than that in \cite{D09}, where a similar statement was proven by cocycle computations under the assumption that $\bar{k}[X]^\times = \bar{k}^\times$.

\begin{lem} \label{comparison H2} With the above notation, one has
$$\partial(\lambda_\sigma ') = 0 \, \, \textup{if and only if} \, \,  \eta'_\sigma = 0 \, .$$
\end{lem} 

\begin{proof} In the following proof, we work over the small \'etale site of $\Spec(k)$. 

Recall that we are given a cocycle $\sigma \in Z^1(k, F)$ as in Lemma \ref{Stoll}: one can associate to $\sigma$ a $\Spec(k)$-torsor $U$ under $F$ with a point $u_0 \in U(\overline{k})$. This torsor $U$ is naturally a homogeneous space of the group $G'$ with geometric stabilizer isomorphic to $T_{\bar{k}}$. 
Section IV.5.1 in \cite{Gi} implies that the element $\eta'_\sigma \in H^2(k, T^\sigma)$ is the class of the $\Spec(k)$-gerbe  $\mathcal{E}_\sigma$ banded by $T^\sigma$ such that for all \'etale schemes $S$ over $\Spec(k)$, the category $\mathcal{E}_\sigma(S)$ is defined as follows: the objects of $\mathcal{E}_\sigma(S)$ are triples $(P,p, \alpha)$ where $P \to S$ is a torsor under $G'$, $p \in P(S_{\bar{k}})$ and $\alpha : P \to U_S$ is a $G'$-equivariant $S$-morphism. The morphisms of $\mathcal{E}_\sigma(S)$ between triples $(P,p,\alpha)$ and $(P',p',\alpha')$ are given by morphisms of torsors $P \to P'$ over $S$ under $G'$ that commute with $\alpha$ and $\alpha'$.

Similarly, one can associate to the morphism $\lambda_\sigma'$ a $\Spec(k)$-gerbe banded by $T^\sigma$ that will be the obstruction for the morphism $\lambda_\sigma'$ to be the extended type of a torsor over $X'$ under $T^\sigma$. The morphism $\lambda_\sigma '$ induces a morphism $\overline{\lambda_\sigma '} : \widehat{T^\sigma_{\bar{k}}} \to KD'(X'_{\bar{k}})$ in $D^b_\et(\bar{k})$. 
By construction, $\overline{\lambda_\sigma '}$ is the extended type of the torsor $Y_0 := W_{\bar{k}} \times_{Z_{\bar{k}}} X'_{\bar{k}}$ over $X'_{\bar{k}}$ under $T^\sigma_{\bar{k}} = T_{\bar{k}}$.

We now define $\mathcal{L}_\sigma$ to be the fibered category defined as follows : for all \'etale schemes $S$ over $\Spec(k)$, the objects of the category $\mathcal{L}_\sigma(S)$ are pairs $(V,\varphi)$, where $V \to X'_S$ is a torsor under $T^\sigma_S$ of extended type $\lambda_V$ compatible with $\lambda_\sigma '$ and $\varphi : V_{\bar{k}} \to Y_0 \times_{\bar{k}} S_{\bar{k}}$ is an isomorphism of torsors over $X' \times_k S_{\bar{k}}$ under $T^\sigma_{S_{\bar{k}}}$. 
Given two such objects $(V, \varphi)$ and $(V', \varphi')$, a morphism between $(V, \varphi)$ and $(V', \varphi')$ in the category $\mathcal{L}_\sigma(S)$ is a pair $(\alpha, t)$, where $\alpha : V \to V'$ is a morphism of torsors over $X'_S$ under $T^\sigma_S$ and $t \in T^\sigma(S_{\bar{k}})$ such that the diagram
\begin{displaymath}
 \xymatrix{
 V_{\bar{k}} \ar[r]^{\overline{\alpha}} \ar[d]_\varphi & V'_{\bar{k}} \ar[d]^{\varphi'} \\
 Y_0 \times_{\bar{k}} S_{\bar{k}} \ar[r]_t & Y_0 \times_{\bar{k}} S_{\bar{k}}
 }
\end{displaymath}
commutes.

One can check that $\mathcal{L}_\sigma$ is a stack for the \'etale topology over $\Spec(k)$, and the fact that this is a gerbe is a consequence of the exact sequence of Proposition 8.1 in \cite{HaSk}
$$H^1(S, T^\sigma) \to H^1(X'_S, T^\sigma) \xrightarrow{\lambda} \Hom_S(\widehat{T^\sigma}, KD'(X'_S)) \xrightarrow{\partial} H^2(S, T^\sigma) $$
(which holds provided that $S$ is integral, regular and noetherian).

The band of this gerbe is the abelian band represented by $T^\sigma$.

In addition, it is clear that $\mathcal{L}_\sigma$ is neutral if and only if $\mathcal{L}_\sigma(k) \neq \emptyset$ if and only if there exists a torsor over $X'$ under $T^\sigma$ of type $\lambda_\sigma '$ if and only if $\partial(\lambda_\sigma ') = 0$.

Let us now construct an equivalence of gerbes between $\mathcal{E}_\sigma$ and $\mathcal{L}_\sigma$.

For all \'etale $\Spec(k)$-schemes $S$, consider the functor
$$m_S : \mathcal{E}_\sigma(S) \to \mathcal{L}_\sigma(S)$$
that maps an object $(P,p, \alpha)$ to the object $(V, \varphi)$, where $V$ is defined to be the contracted product $V := (P \times_S^{G'} W_S) \times_{Z^\sigma_S} X'_S$ and $\varphi : V_{\bar{k}} \to Y_0 \times_{\bar{k}} S_{\bar{k}} = (W_{\bar{k}} \times_{Z_{\bar{k}}} X_{\bar{k}}') \times_{\bar{k}} S_{\bar{k}}$ is induced by the point $p \in P(S_{\bar{k}})$. Indeed, by construction, we have a natural map $P \times_S^{G'} W_S \to U_S \times_S^F Z_S = Z^\sigma_S$, and a simple computation proves that this map is a torsor under $T^\sigma$ of extended type compatible with $\lambda_W^\sigma$.

By definition, the functor $m_S$ sends a morphism $\varphi : (P,p, \alpha) \to (P', p', \alpha')$ to the morphism $(\widetilde{\varphi}, t_0)$ such that $\widetilde{\varphi} : (P \times_S^{G'} W_S) \times_{Z^\sigma_S} X'_S \to (P' \times_S^{G'} W_S) \times_{Z^\sigma_S} X'_S$ is the morphism induced by the morphism of torsors $\varphi : P \to P'$, and $t_0 \in T^{\sigma}(S_{\overline{k}})$ is the element such that $p' = t_0 \cdot \varphi(p)$ as $S_{\overline{k}}$-points in $(P' \times_S^{G'} W_S) \times_{Z^\sigma_S} X'_S$.

Finally, one checks that the collection of functors $m_S$ defines a morphism of gerbes $m : \mathcal{E}_\sigma \to \mathcal{L}_\sigma$ banded by the identity of $T^\sigma$, which implies that $\eta'_\sigma := [\mathcal{E}_\sigma] = [\mathcal{L}_\sigma] \in H^2(k, T^\sigma)$.

Therefore, $\eta'_\sigma = 0$ if and only if $\mathcal{E}_\sigma(k) \neq \emptyset$ if and only if $\mathcal{L}_\sigma(k) \neq \emptyset$ if and only if $\partial(\lambda_\sigma ') = 0$.
\end{proof}

The immediate consequence of Lemma \ref{comparison H2} is the following result which extends Proposition 5 in \cite{D09} to open varieties. 

\begin{prop}
\label{prop}
Let $X$ be a smooth geometrically integral $k$-variety. Let $(P_v) \in X({\bf A}_k)^{\textup{\'et, Br}}$ and let $Y \rightarrow X$ be a torsor under a linear $k$-group $G$.
Let
$$1 \rightarrow H \rightarrow G \rightarrow F \rightarrow 1$$
be an exact sequence of linear $k$-groups, where $H$ is connected and $F$ finite. Let $Z \rightarrow X$ be the push-forward of $Y \rightarrow X$ by the morphism $G \rightarrow F$, which is a torsor under $F$. Let $\sigma \in Z^1(k, F)$ be a $1$-cocycle given by Lemma \ref{Stoll} applied to the torsor $Z \rightarrow X$ and the point $(P_v)$.

Then the cocycle $\sigma \in Z^1(k, F)$ lifts to a cocycle $\tau \in Z^1(k, G)$.
\end{prop}

\begin{proof} As mentionned above, Construction (5.1) in \cite{FSS} (see also (7.7) in \cite{B}) gives a class $\eta_\sigma$ of $H^2(k, \kappa_\sigma)$ such that $\sigma$ can be lifted to $Z^1(k,G)$ if and only if $\eta_\sigma$ is neutral, where $\kappa_\sigma$ is a $k$-kernel on $H_{\bar{k}}$. By (6.1.2) of \cite{B} and Lemma 6 in \cite{D09}, there is a canonical map $H^2(k, \kappa_\sigma) \to H^2(k, T^\sigma)$ such that the class $\eta_\sigma$ is neutral if and only if its image $\eta'_\sigma \in H^2(k, T^\sigma)$ is zero. By Lemma \ref{comparison H2}, one only needs to show that $\partial(\lambda_\sigma ') = 0$ where $\lambda_\sigma '=\psi^* \circ \lambda_W^\sigma$, with $KD'(Z^\sigma)\xrightarrow{\psi^*} KD'(X')$ given by Lemma \ref{Stoll} and $\lambda_W^\sigma$ defined by Lemma \ref{twist}. 

By Lemma \ref{Stoll}, we know that $X'({\bf A}_k)^\Br \neq \emptyset$. Therefore the map $\lambda$ in the exact sequence (see Proposition 8.1 in \cite{HaSk})
$$H^1(X', T^\sigma) \xrightarrow{\lambda} \Hom_k(\widehat{T^\sigma}, KD'(X')) \xrightarrow{\partial} H^2(k, T^\sigma) $$
is surjective by Corollary 8.17 in \cite{HaSk}. Hence the map $\partial$ is the zero map and $\partial(\lambda_\sigma ') = 0$, which concludes the proof. \end{proof}

\begin{rem} The proof of Proposition \ref{prop} also gives the following result: 

Let $X$ be a smooth geometrically integral $k$-variety and let $Y \rightarrow X$ be a torsor under a linear algebraic $k$-group $G$. Let
$$1 \rightarrow H \rightarrow G \rightarrow F \rightarrow 1$$
be an exact sequence of linear $k$-groups, where $H$ is connected and $F$ finite. Let $Z \rightarrow X$ be the push-forward of $Y \rightarrow X$ by the morphism $G \rightarrow F$.

If $\sigma\in H^1(k, F)$ satisfies $Z^{\sigma}({\bf A}_k)^{\Br_1(Z^\sigma)}\neq \emptyset$, then $\sigma$ can be lifted to $H^1(k, G)$.
\end{rem}

One can now prove the main result of this section:

\begin{thm}\label{main-last} If $X$ is a smooth and geometrically integral variety over a number field $k$, then
$$ X({\bf A}_k)^{\textup{\'et}, \Br} \subseteq  X({\bf A}_k)^{\rm{desc}} \, . $$
\end{thm}

\begin{proof} 
Since the statement 2 of Theorem 2 in \cite{H} (which we apply to $X'$) holds for any geometrically integral variety (without any assumption on $\bar{k}[X']^\times$), the proof of this theorem using Proposition \ref{prop} is exactly the same as the proof of Theorem 1 using Proposition 5 in \cite{D09} (see in particular \cite{D09}, p. 244-245). \end{proof}


\bigskip

\noindent{\bf Acknowledgements.} We would like to thank Jean-Louis Colliot-Th\'el\`ene for several comments for the early version of this paper. We would also like to thank the referee for pointing out a mistake in the previous version. 
 The first named author acknowledges the support of the French Agence Nationale de la Recherche (ANR)
 under reference ANR-12-BL01-0005, the second named author acknowledges the support of the French Agence Nationale de la Recherche (ANR)
 under references ANR-12-BL01-0005 and ANR-15-CE40-0002-01, and the third named author acknowledges the support of  NSFC grant no.11471219 and 11631009.


\end{document}